\documentstyle[a4wide]{article}
\newtheorem{lem}{Lemma}
\newtheorem{thm}{Theorem.}
\newcommand{\dfrac}[2]{{\displaystyle \frac{#1}{#2}}}
\newcommand{\infrac}[2]{{ {#1}/{#2}}}

\newcommand{\samplex}{$X_{1},\ldots,X_{m}$}
\newcommand{\sampley}{$Y_{1},\ldots,Y_{n}$}

\newlength{\longueur}
\newcommand{\B}[1]{B_{#1}^{o}}
\def\k#1{\kern#1em}
\def\Ib#1{{I\kern-.25em#1}}
\def\Ibb#1{{I\kern-.23em#1}}

\newcommand{\DLZZ}{\settowidth{\longueur}{Z}Z\hspace{-0.9\longueur}Z}

\begin{document}

\begin{titlepage} 
\title{Dominance refinements of the
Smirnov two-sample test}
\author{A. Di Bucchianico$^{*}$  \\
Department of Mathematics and Computing Science\\ 
Eindhoven University of Technology \\
P. O. Box 513\\
5600 MB Eindhoven, The Netherlands\\
{\tt sandro@win.tue.nl}
\and 
D.\ E.\ Loeb\thanks{Authors supported by 
NATO CRG 930554. Second author partially supported by URA CNRS 1304, EC grant
CHRX-CT93-0400, and the PRC Maths-Info.}\\
LaBRI\\ 
Universit\'e de Bordeaux I\\
33405 Talence, France\\
{\tt loeb@labri.u-bordeaux.fr}\\
URL 
{\tt http://www.labri.u-bordeaux.fr/$\sim$loeb}} 
\date{ }
\end{titlepage} 
\maketitle

\begin{abstract} 
\noindent
We prove the following conjecture of Narayana: there are no nontrivial 
dominance refinements of the Smirnov two-sample test if and only if the two sample 
sizes are relatively prime. 
We also count the number of natural significance levels of the Smirnov two-sample test 
in terms of the sample sizes and relate this to the Narayana conjecture. In particular, 
Smirnov tests 
with relatively prime sample sizes turn out to have many more natural significance 
levels than do Smirnov tests whose sample sizes are not relatively prime (for example, 
equal sample sizes). 
\end{abstract}


\section{Introduction}
Let \samplex\ and \sampley\ be independent random samples from continuous distribution 
functions $F$ and $G$, respectively. In order to test
nonparametrically whether  $X_{1}$ is 
stochastically smaller than $Y_{1}$, one often uses the {\bf Smirnov 
statistic} 
\begin{equation}\label{dmnplus}
D_{mn}^{+} = \sup_{t}\, (F_{m}(t) - G_{n}(t) ),
\end{equation}
where $F_{m}$ and $G_{n}$ are the empirical distribution functions of \samplex\ 
and \sampley\ respectively. 

Narayana (1971, pp. 43 ff.) constructs explicit dominance refinements of
(upper-tailed) Smirnov 
two-sample tests  with equal sample sizes and calculates their power against 
Lehmann alternatives $G = F^k (k
> 0)$. In
the ranges considered by Narayana ($3 \leq m = n \leq 10$), he shows that the powers
of these dominance refinements are uniformly greater than the Smirnov test. It is thus 
of practical importance to know when dominance refinements exist. 

Narayana (1975) stated without proof that dominance refinements of the
Smirnov two-sample test exist if and only if $\gcd (m,n) > 1$. This claim
was restated as a conjecture in Narayana (1979, Exercise~9, p.~30). The
purpose of this paper is to prove this conjecture. We will show 
that this conjecture is closely related to the number of natural significance
levels of the Smirnov two-sample test.

\section{Dominance refinements}

In this section, we explain what dominance refinements are.

A convenient way to study the distribution
of $D_{mn}^{+}$ is  
the so-called Gnedenko path. The {\bf Gnedenko path} $\omega$ of the samples
\samplex\ and \sampley\ is 
defined as follows: $\omega$ is a path from $(0,0)$ to $(m,n)$ with unit
steps $\omega_{i}$ to the east or north. If the $i$th value of the
ordered combined sample comes from \samplex, then $\omega_{i}$ is a
step east; otherwise, it is a step north.  Since we assume that $F$
and $G$ are continuous, 
the probability of a tie ({\em i.e.}, $X_i=Y_j$) is zero. Hence, $\omega$
is  almost surely well-defined. It is easy to see that under $H_{0}: F = G$, all 
paths from 
$(0,0)$ to $(m,n)$ are equiprobable, {\em i.e.}, ${\bf P}(w) = 1/{m+n \choose n}$ 
for all paths $w$. Now, 
$$  m n D_{mn}^{+} \leq r  $$
if and only if all vertices $(x,y)$ of the Gnedenko path satisfy 
$$nx - my \geq r$$
In other
words, $ m n D_{mn}^{+} \geq r$ if and only if $\omega$ crosses
below 
the line $nx - my = r$. A convenient way to describe a 
path $\omega$ is to represent it by a $n$-tuple $\langle
t_{1},\ldots,t_{n} \rangle$, where 
$t_{i}$ is the horizontal distance from $(m,n-i)$ to $\omega$
(see Figure~1).
\begin{figure}
\label{fig1}
\begin{center}
\begin{picture}(5,4)
\put(0,0){\circle*{.1}} 
\put(0,0){\line(0,1){1}}
\put(0,1){\circle*{.1}} 
\put(0,1){\line(1,0){2}}
\put(1,1){\circle*{.1}} 
\put(2,1){\circle*{.1}} 
\put(2,1){\line(0,1){2}}
\put(2,2){\circle*{.1}} 
\put(2,3){\circle*{.1}} 
\put(3,3){\circle*{.1}} 
\put(4,3){\circle*{.1}} 
\put(2,3){\line(1,0){2}}
\thicklines
\put(4,2){\vector(-1,0){2}}
\put(4.3,2){$t_{1}=2$} 
\put(4,1){\vector(-1,0){2}}
\put(4.3,1){$t_{2}=2$} 
\put(4,0){\vector(-1,0){4}}
\put(4.3,0){$t_{3}=4$} 
\end{picture}
\end{center}
\caption{Representation of a Path}
\end{figure}
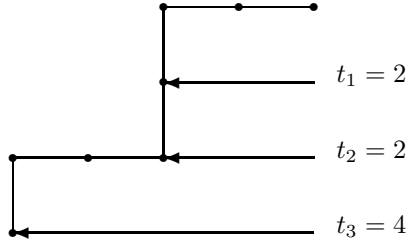
The path $\langle s_{1},\ldots,s_{n} \rangle $ is said to {\bf dominate} 
$\langle t_{1},\ldots,t_{n} \rangle  $ if 
$s_{i} \geq t_{i}$ for $i = 1,\ldots,n$. 

Let the $\lceil x \rceil$ be the {\bf ceiling} of $x$, {\em i.e.}, 
the smallest integer larger 
than or equal to $x$ and let $\lfloor x \rfloor$ be the {\bf floor} of
$x$, {\em i.e.},   
the largest integer 
not exceeding $x$. The {\bf $r$-profile} is the path  
$\langle \max(0,t_{1}),\ldots,\max (0,t_{n}) \rangle  $, where $t_i$ is the
ceiling of the  
horizontal distance from $(x,m-i)$ to the line $nx - my = r$). Clearly, the
$r$-profile is the minimal path that lies above 
(possibly touching) the line  $nx - my = r$.

Thus, we may cast the (upper-tailed) Smirnov two-sample test completely in terms of 
Gnedenko paths as 
follows: $ m n D_{mn}^{+} \leq r$ if and only if the Gnedenko path dominates the 
$r$-profile. Thus, the Smirnov two-sample test is completely characterized by its
$r$-profiles ({\em i.e.}, we regard the test as a set of critical regions, indexed by
its natural significance levels, cf. Gibbons (1992, p.~23)). This formulation shows 
that we attain more significance levels 
if we can insert intermediate paths between consecutive $r$-profiles of the  Smirnov 
two-sample test (see Narayana (1979, Chapter~2)). A set of paths totally ordered
by dominance is said to be a 
{\bf dominance refinement} of any set of paths included in it. Note that under this
definition, we consider a test to be a 
dominance refinement of itself, called the  {\bf trivial dominance refinement}. A set 
of paths is {\bf saturated} if it has no nontrivial
dominance refinement.  

Of course, there exist other refinements of  the Smirnov
test. Each partition of the set of paths with a common value $r$ of the statistic
$D_{mn}^{+}$ ({\em i.e.}, all paths that touch but do not cross the line $nx - my = r$)
yields a refinement of the Smirnov test. For example, we can divide the paths that 
touch but
do not cross the line $nx - my = r$ according to the number of times that they touch 
the line $nx - my = r$. Dominance refinements partition the set of paths with a common
value of $D_{mn}^{+}$ into dominance regions, {\em i.e.}, collections
of paths that dominate 
a given path. An advantage of dominance refinements is that they can be
described very efficiently by simply listing the critical paths. Hence, the
refined test can
be performed graphically. 

Another reason for considering dominance refinements (or the
notion of dominance itself) is the following relation with most powerful 
rank (MPR) tests.  
If $F$ and $G$ have densities $f$ and $g$ respectively, and the
likelihood ratio $f/g$ is 
increasing (as is the case for the Lehmann alternatives $H_{a} : G = F^{k}$, $k >
0$), then $s$
dominates $t$ implies ${\bf P}(t | G = F^{k}) \geq {\bf P}(s | G = F^{k})$ (see
Savage (1956)). Hence, if a path $s$ belongs to the critical region of an MPR test, 
then all paths dominated by $s$ must also belong to this critical region. Thus, an 
MPR test at a fixed
significance level is a
dominance test in the terminology of  Narayana (1979, Chapter~3,
p.~35). Conversely, dominance tests are good approximations for MPR
tests (see Narayana (1979, Chapter~3, pp. 44-45)).

\begin{center}
\begin{tabular}{l@{\hspace{1in}}l}
\begin{picture}(5,3)
\put(4.3,2.0){\footnotesize $r=0$}
\put(4.3,1.5){\footnotesize $r=2$}
\put(4.3,1.0){\footnotesize $r=4$}
\put(4.3,0.5){\footnotesize $r=6$}

\thinlines
\put(0,0){\line(2,1){4}} 
\put(1,0){\line(2,1){3}} 
\put(2,0){\line(2,1){2}} 
\put(3,0){\line(2,1){1}} 

\put(0,0){\line(1,0){1}}
\put(1,0){\line(0,1){1}}
\put(2,1){\line(1,0){1}}
\put(3,1){\line(0,1){1}}
\newsavebox{\mydots}
\savebox{\mydots}{\multiput(0,0)(1,0){5}{\circle*{.1}}}
\multiput(0,0)(0,1){3}{\usebox{\mydots}}
\thicklines
\put(0,0){\line(0,1){1}}
\put(0,1){\line(1,0){2}}
\put(2,1){\line(0,1){1}}
\put(2,2){\line(1,0){2}}
\end{picture}
&
\begin{picture}(6,4)
\put(5.3,3.0){\scriptsize $r=0$}
\put(5.3,2.8){\scriptsize $r=1$}
\put(5.3,2.6){\scriptsize $r=2$}
\put(5.3,2.4){\scriptsize $r=3$}
\put(5.3,2.2){\scriptsize $r=4$}
\put(5.3,2.0){\scriptsize $r=5$}
\put(5.3,1.8){\scriptsize $r=6$}
\put(5.3,1.6){\scriptsize $r=7$}
\put(5.3,1.2){\scriptsize $r=9$}
\put(5.3,1.0){\scriptsize $r=10$}
\put(5.3,0.6){\scriptsize $r=12$}

\put(0,0){\line(5,3){5}}         
\put(0.333,0){\line(5,3){4.667}} 
\put(0.667,0){\line(5,3){4.333}} 
\put(1,0){\line(5,3){4}}         
\put(1.333,0){\line(5,3){3.667}} 
\put(1.667,0){\line(5,3){3.333}} 
\put(2,0){\line(5,3){3}}         
\put(2.333,0){\line(5,3){2.667}} 
\put(3,0){\line(5,3){2}}         
\put(3.333,0){\line(5,3){1.667}} 
\put(4,0){\line(5,3){1}}         

\put(1,1){\line(1,0){1}}
\put(2,1){\line(0,1){1}}

\newsavebox{\toy}
\savebox{\toy}{\multiput(0,0)(1,0){6}{\circle*{.1}}}
\multiput(0,0)(0,1){4}{\usebox{\toy}}
\thicklines
\put(0,0){\line(0,1){1}}
\put(0,1){\line(1,0){1}}
\put(1,1){\line(0,1){1}}
\put(1,2){\line(1,0){2}}
\put(3,2){\line(0,1){1}}
\put(3,3){\line(1,0){2}}
\end{picture}\\[3mm]
Figure 2: $m=4$ and $n=2$ & Figure 3:  $m=5$ and $n=3$.
\end{tabular}
\end{center}
Let us look at two examples in order to get a feeling for the Narayana
conjecture. 
\begin{itemize}
\item {\bf (Figure 2, $m=4$, $n=2$)} The 0-profile and
the 1-profile coincide and are equal to the path  
$\langle 2,4 \rangle  $
whereas the 2-profile is the path $\langle 1,3 \rangle  $.
Thus, we see that there are two intermediate paths between the 1-profile and
the 2-profile: $\langle 1,4 \rangle  $
and $\langle 2,3 \rangle  $. Inserting either of these paths, we
obtain a refinement of the Smirnov test. 
Note that the 2-profile differs from the 1-profile by the possibility to go
through the points $(1,0)$ and $(3,1)$, which both lie on the line $2x - 4y = 2$. 
\item {\bf (Figure 3, $m=5$, $n=3$)} The 0-profile is the path
$\langle 2,4,5 \rangle $, and the  
1-profile is the path $\langle 2,3,5 \rangle  $. Thus, there is no
intermediate path between the  
profiles; this is also true for the other pairs of consecutive
profiles. In other words, there is no 
refinement. Note that there is only one lattice point on each line of the form 
$3x - 5y = r$ and that no profiles coincide.   
\end{itemize}

\section{Main Results}

The examples above 
indicate that the existence of dominance refinements
depends on the  
number of lattice points on lines of the form $nx - my = r$.
In the following lemmas, we enumerate these points. These lemmas are
used in our  proof of the 
Narayana conjecture. 

For unexplained notions of number theory and
combinatorics,  we refer the reader to Stark (1970, Chapters~2 and 3) 
and Berge (1971).

\begin{lem}\label{L1}
Let $m$, $n$, and $r$ be positive integers, and let $d=\gcd(m,n)$. 

{\bf 1.} The Diophantine equation
\begin{equation}\label{unconstrained}
n x - m y = r
\end{equation}
has integer solutions if and only if $d$ divides  $r$.

{\bf 2.} If $(x,y)$ is a solution of (\ref{unconstrained}), then
the integer solutions of  
(\ref{unconstrained}) are exactly 
\begin{equation}\label{steps}
\left\{\left(x + t \dfrac{m}{d},y + t
\dfrac{n}{d}\right): t \in \DLZZ \right\}. 
\end{equation}
\end{lem}

{\em Proof}: {\bf 1.} If $r = d$, then by Euclid's Lemma there exist
integer 
solutions $(x,y)$ of $n x - m y = r$. Obviously, this also holds if $r$ is a multiple 
of $d$.
Conversely, if there exists an integer solution $(x,y)$ of $n x - m y = r$, then $r$ 
is a multiple of $d$, since $d$ divides both $m$ and $n$.

{\bf 2.} If $t$, $x$, and $y$ are integers such that $n x - m y = r$,
then $x' := x + t 
\infrac{m}{d}$ and  $y' := y + t \infrac{n}{d}$  satisfy  $n x' - m y' = r$. Conversely,
if  $n x - m y = r$ and  $n x' - m y' = r$, then subtraction yields 
$n(x - x') = m (y' - y)$. Cancelling the common factor $d$ and using the uniqueness of
prime factorizations, we see that there exists an integer $t$ such that $x - x' = t
\infrac{m}{d}$ and  $y' - y = t \infrac{n}{d}$. \hfill $\Box$

\begin{lem}\label{previous}
Let $m$ and  $n$ be positive integers with greatest common divisor
$d$ and $m\geq n$. Let $r$ be a 
nonnegative number. 
The Diophantine equation
\begin{equation}\label{star}
n x - m y = r\ \ \ ( 0 \leq x \leq m\ \mbox{and}\ 0 \leq y \leq n)
\end{equation}
has integer solutions only if $d$ divides $r$ and $0\leq r\leq nm.$
In that case,  the number of
solutions to the Diophantine equation (\ref{star}) 
is given by 
\begin{equation}\label{sols}
\left.
\begin{array}{lcr}
& &\alpha_r :=  d + 1 - \left\lceil \dfrac{p +a}{m/d} \right\rceil,\\
where & &\\
& & p = \left\lfloor \dfrac{r/d}{n/d} \right\rfloor = 
\left\lfloor \dfrac{r}{n} \right\rfloor\\
and & &\\
& & a =  \dfrac{r-n p}{d}\ (n/d)^{-1}\ \bmod\ m/d.
\end{array}
\right\}
\end{equation}
Here, 
$(n/d)^{-1}$ denotes the inverse of $n/d$ in the ring 
 ${\DLZZ}_{m/d} =
\left\{ 0,1,\ldots,m/d -1 \right\} $.
\end{lem}

{\em Proof}: The first assertion is an immediate consequence of Lemma
\ref{L1} and the fact that $(m,0)$ lies on the line $nx - my = nm$.

Now, 
define $m' = m/d$, $n' = n/d$, and $r' = r/d$. 
By the remainder theorem, there exist nonnegative integers $p$ and $q$ such that 
$r' = n' p + q$ with $0 \leq q < n'$. Since $m'$ and $n'$ are relatively prime, 
it follows from Lemma~1 that there exist integers 
$a\ (0 \leq a \leq m')$ and $b\ (0 \leq b < n')$ such that $q = a n' - b
m'$. Thus, $r' = n'(p + a) - m' b$. Hence, $(p+a,b)$ is a solution to 
(\ref{unconstrained}). We claim that 
$(p+a,b)$ is the solution $(x,y)$ $(x \geq 0)$ of
(\ref{unconstrained}) with minimal  
nonnegative $y$. This claim follows from the above, since the next
smaller solution of 
(\ref{unconstrained})
has $y' = b - n' < 0$. Thus, the solutions of (\ref{unconstrained})
with $y \geq 0$  
are exactly $\left\{ (p + a + t m', b + t n') \right\}_{t \geq 0}$.

Now, 
impose $x \leq m$ on this solution set. This leads to $t \leq (m - p - a)/m'$.
Therefore, the solutions to (\ref{star}) are 
$\left\{ (p + a + t m', b + t n') : 0 \leq t \leq T \right\}$ 
 where $T = \lfloor \frac{m -
p -a}{m'} \rfloor = d - \lceil \frac{p+a}{m'} \rceil$. 

Note that if $ p+a > m$, then $T =-1$ and there are no solutions at all.

The value of $p$ follows from  $ n' p =  r' - q$ ($0 \leq q < n'$). It  
 follows from $r' = n' (p + a) - m' b$ that $n' a \equiv  r' - n' p
\bmod m'$. Since $m'$ 
and $n'$ are relatively prime, $n'$ has an inverse in  ${\DLZZ}_{m'}$.
Thus, $a \equiv (r' - n' p)\, (n')^{-1} \bmod m'$. Equivality follows
since $a<n'\leq m'.$
\hfill$\Box$\bigskip\bigskip 

\newpage
\noindent
{\bf Examples.}\begin{itemize}\nopagebreak[4]
\item {\bf (Figure 2, $m=4$, $n=2$)} 
If $r=1$, then there are no solutions, since $d$ does 
not divide $r$. 
If $r = 2$, then $p = \lfloor 1/2 \rfloor = 0$ and 
$n/d = 1$, hence $(n/d)^{-1} \equiv 1 \bmod 2$ and $ a = 0$. Thus,
the number  
of solutions equals $2 + 1 - \lceil \frac{1+0}{2} \rceil = 2 + 1 - 1 = 2$.
\item {\bf (Figure 3, $m=5$, $n=3$)} If $r = 1$, then $p = \lfloor 1/3
\rfloor = 0$ and $(n/d)^{-1} \equiv  3^{-1} \equiv 2 \bmod 5$, 
since $6 \equiv 1 \bmod 5$. 
Thus, 
 $a = 2$ and the number of solutions equals $1 + 1 - \lceil
\frac{0+2}{5} \rceil = 1 + 1 -0 = 1$. If $r = 14$, then $ p = \lfloor
14/3 \rfloor = 
4$, $ a = 2$, and the number of solutions equals 
$ 1 + 1 - \lceil 6/5 \rceil = 1+1-2 = 0$.
\end{itemize}

\begin{lem}\label{last}
Let $m$ and $n$ be positive integers with greatest common divisor $d$. 
Let $s_{k}$ be the number of positive divisors $r$ of $d$  
such that $\alpha_{r}=k$ (in Lemma \ref{previous}), in other words, 
such that the Diophantine equation \ref{star} has
precisely $k$ 
integer solutions. Then
\begin{itemize}
\item[\bf 1.] $s_{k}=0$ for $k>d$,
\item[\bf 2.] $s_{k}={nm}/{d^2}$ for $0<k<d$,
\item[\bf 3.] $s_{0}=(nm-(n+m)d+d^{2})/2d^{2}$, and 
\item[\bf 4.] $s_{d}=(nm+(n+m)d-d^{2})/2d^{2}$.
\end{itemize}
\end{lem} 

{\em Proof}: We use the notation of Lemma~2 and its proof. Without
loss of generality, $m\geq n$. 

{\bf 1.} This follows from the second part of Lemma~1.

{\bf 2.} By Lemma~2, $n x - m y =r$ has $k$ ($0 < k < d$) solutions with the 
constraints 
$ 0 \leq x \leq m$ and $0 \leq y \leq n$ if $\lceil \frac{p + a}{m'} \rceil = 
d + 1 -k$. In other words,
if $m - m' k < p + a \leq m - m' (k-1)$. Each $r'$ corresponds uniquely to a pair
$(p,q)$ with either $0 \leq q < n'$ and $0 \leq p < m$ or $q = 0$ and $p = m$.
First choose $q$, which can be done in $n'$ ways. This fixes $a$. Hence, $p$ must obey
$0 \leq m - m'k - a < p \leq  m - m' (k-1) - a \leq m -a \leq m$, of which there are
exactly $m'$ solutions. 

{\bf 3.} By Lemma~1, we must have $r > nm - nm/d$, or equivalently 
$r' > nm/d - nm/d^{2}$. Each
line $n'x - m'y = r'$ has 0 or 1 integer solutions with $ 0 \leq x \leq m$ and 
$0 \leq y \leq n$. Now, consider all lattice points $(x,y)$ with 
$n x - m y > n m - nm/d$,
$x \leq m$, and $y \geq 0$, {\em i.e.}, 
 the triangle spanned by the points $(m-m',0)$, 
$(m,0)$, and
$(m,n')$ with the points  $(m-m',0)$ and $(m,n')$ excluded. It is easy to see that this 
triangle
contains $\Delta := \infrac{(m'+1) (n'+1) -2)}{2}$ points. By the second part of 
Lemma~1, each such point 
corresponds uniquely to a line $nx - my = r$ with exactly one lattice point in the 
region 
$ 0 \leq x \leq m$ and $0 \leq y \leq n$. Hence, the other $n m /d^2 - \Delta$ lines 
contain no point.

{\bf 4.} This follows by elimination, 
since there are $nm/d$ values of $r$ under consideration.\hfill $\Box$

Define the {\bf  Bell-Stirling number} $\B{n}$ to be the number of
{\bf ordered partitions} of $\{1,2,\ldots,n \}$ into disjoint nonempty
sets $S_{1}\cup \cdots \cup
S_{k}=\{1,2,\ldots,n \}$. 
For example, $\B{0}=1,$ $\B{1}=1,$ $\B{2}=3,$ $\B{3}=13,$ $\B{4}=75,$ 
$\B{5}=539$.
The  Bell-Stirling number can also be defined implicity by its
generating function (Motzkin (1971, p. 171))
$$ \sum _{n=0}^{\infty } \B{n}t^{n}/n! = \frac{1}{2-\exp(t)}.$$

Now, we prove the Narayana conjecture by counting 
dominance  refinements and saturated dominance refinements. 
\begin{thm}
Let $m$ and $n$ be positive integers with greatest common divisor $d$.
\begin{enumerate}
\item[\bf 1.] The Smirnov upper-tailed two-sample test 
with sample sizes $m$ and $n$ has exactly $(d^2+ nm(2d - 1) +d (n+m))/(2d^{2})$
natural levels.
\item[\bf 2.] All saturated dominance refinements of the 
Smirnov upper-tailed two-sample test with sample sizes $m$ and $n$ have
exactly $((n+1) (m+1)-(d+1))/2 + 1$ natural levels.
\item[\bf 3.] {\bf (Narayana's Conjecture)}  
The Smirnov upper-tailed two-sample test is saturated if and only if 
the sample sizes are
relatively prime.
\item[\bf 4.] The number of dominance refinements of 
the Smirnov upper-tailed two-sample test with samples size $m$ and $n$
(including the trival one) is given by the product
$$(\B{d})^{(nm+(n+m)d-d^{2})/( 2d^{2})} 
\prod_{k=1}^{d-1} (\B{k})^{nm/d^{2}}.$$  
\item[\bf 5.] The number of saturated dominance refinements of 
the Smirnov upper-tailed two-sample test with sample sizes $m$ and $n$
is given by the product
$${d}!^{(nm+(n+m)d-d^{2})/(2 d^{2})} 
\prod_{k=1}^{d-1} {k}!^{nm/d^{2}}.$$  
\end{enumerate}
\end{thm}
{\em Proof}: {\bf 1.} We use the representation of the Smirnov test as
a set of $r$-profiles. As noted above, the $r$-profile is different
from the $(r-1)$-profile exactly when the Diophantine equation
(\ref{star}) has at least one solution $(x,y)$. The number of levels
thus corresponds to the number of $r$ such that $\alpha_{r}>0$ in
Lemma \ref{previous} (including $r = 0$). We are
thus led to sum $1 + s_{1}+s_{2}+\cdots +s_{d}$ in Lemma \ref{last}, which equals 
$1 + (d-1) \infrac{nm}{d^2} + \infrac{(n m+(n+m) d - d^2)}{(2 d^2)} = 
\infrac{(d^2+ nm(2d - 1) +d (n+m))}{(2d^{2})}$.

{\bf 2.} Let $(P_{0},P_{1},\ldots,P_{k})$ be a saturated set of paths.
The region strictly between two consecutive paths $P_{i}$ and $P_{i+1}$ is
exactly a unit square 
$\{(x,y) \mid (x_{i}-1 \leq x_{i} \mbox{and} y_{i} \leq y_{i}+1\}$ 
 as in Figure~2, otherwise additional paths
could be inserted between $P_{i}$ and $P_{i+1}$. 
All points $(x,y)$ with $0\leq x\leq m$, $0\leq y\leq n$ and
$nx-my > 0$ determine exactly one such unit square, since they lie
under the 0-profile. There are 
$((n+1)(m+1) - (d + 1))/2$ 
such points, since there are
$d+1$ lattice points on $nx - my = 0$. Thus,
the path is of length $(n+1)(m+1)- (d + 1))/2-1$ 
and consists of 
 $(n+1)(m+1)- (d + 1))/2$ 
natural 
significance 
levels.

{\bf 3.} If $m$ and $n$ are relatively prime, then substitution of $d=1$ into 
parts 1) and 2) yields that the Smirnov upper-tailed two-sample test has as many
natural significance levels as its saturated dominated refinement. Hence, the Smirnov
test is saturated. Conversely, if $m$ and $n$ are not relatively prime, then by 
Lemma~\ref{last} there is at least one line $nx - my = r$ with two lattice
points. Hence, the Smirnov upper-tailed two-sample test admits a non-trivial dominance
refinement.

{\bf 4.} Consider an $r$-profile $P$ with $\alpha_{r}=k$, and the
$(r-1)$-profile $P'$. 
The area between $P$ and $P'$ consists of $k$ unit squares as in Figure~3.
New profiles can be ``inserted'' between $P$ and $P'$ by passing the
path across the $k$ unit squares in several steps instead of all at
once. Let $S_{1}\cup \cdots \cup S_{k}$ be an ordered partition of the
$k$ unit squares. Define the intermediary path $P_{i}$ to pass
under the cells $S_{1}\cup \cdots \cup S_{i}$, over the cells
$S_{i+1}\cup \cdots \cup S_{k}$, and follow paths $P$ and $P'$ where
they agree. Each such ordered partition determines a dominance
refinement of this step in the list of paths which make up the Smirnov test.
There are $\B{k}$ ordered partitions.

To consider the combinations of refinements of each step of the Smirnov test, 
we take the product over all $r$.

{\bf 5.} The reasoning is similar to that in part 4 except we must
only consider ordered partitions into unit blocks. There are $k!$ such
ordered partitions of a $k$ element set.\hfill $\Box $

Table~\ref{natu} shows the number of natural significance levels for the
upper-tailed Smirnov two-sample test for various sample sizes. Note especially the
low number of levels for Smirnov tests with equal sample sizes ({\em e.g.}, 
compare the number of levels for $m=n=10$ with those for $m=10$ and $n=9$). 
Table~\ref{satur} shows the number of natural significance levels for saturated 
dominance refinements of the
upper-tailed Smirnov two-sample test for various sample sizes. Note that unlike
in Table~\ref{natu}, there are no big differences in Table~\ref{satur}.

\begin{table}
\begin{center}
\begin{tabular}{|c|rrrrrrrr|} \hline
$ m\backslash n$ &  3 &  4 &  5 &  6  &  7  &   8 &   9 &  10 \\ \hline
3     &  4 &  10 & 12 &  7 &  16 &  18 &  10 &  22 \\
4     & 10 &   5 & 15 & 12 &  20 &   9 &  25 &  19 \\
5     & 12 &  15 &  6 & 21 &  24 &  27 &  30 &  11 \\
6     &  7  & 12 &  21&  7 &  28 &  22 &  18 &  27 \\
7     & 16  & 20 &  24& 28 &   8 &  36 &  40 &  44 \\
8     & 18  &  9 &  27& 22 &  36 &   9 &  45 &  35 \\
9     & 10  & 25 &  30& 18 &  40 &  45 &  10 &  55 \\
10    & 22  & 19 &  11& 27 &  44 &  35 &  55 &  11 \\ \hline
\end{tabular}
\caption{Number of natural significance levels 
of the upper-tailed Smirov two-sample test.}\label{natu}
\end{center}
\end{table}
\begin{table}
\begin{center}
\begin{tabular}{|c|rrrrrrrr|} \hline
$ m\backslash n$ &  3 &  4 &  5 &  6  &  7  &   8 &   9 &  10 \\ \hline
3     &  7 &  10 & 12 &  13 &  16 &  18 &  19 &  22 \\
4     & 10 &  11 & 15 &  17 &  20 &  21 &  25 &  27 \\
5     & 12 &  15 & 16 &  21 &  24 &  27 &  30 &  31 \\
6     & 13  & 17 &  21&   22&   28&   31&   34&   38 \\
7     & 16  & 20 &  24&   28&   29&   36&   40&   44 \\
8     & 18  & 21 &  27&   31&   36&   37&   45&   49 \\
9     & 19  & 25 &  30&   34&   40&   45&   46&   55 \\
10    & 22  & 27 &  31&   38&   44&   49&   55&   56 \\ \hline
\end{tabular}
\caption{Number of natural significance levels of saturated dominance refinements
of the upper-tailed Smirov two-sample test.}\label{satur}
\end{center}
\end{table}

If we can list the profiles of an upper-tailed Smirnov two-sample test, then using
Kreweras' theorem (Narayana
(1971, p. 21) we can also calculate the levels rather  
than the total number of levels. For example, in the case of equal sample sizes $m=n$,
the profiles are the paths 
$\langle 1,2,\ldots, m\rangle$, $\langle 0,1,2,\ldots, m-1\rangle$,$\ldots$, 
$\langle 0,0,\ldots, 0\rangle$. In order to calculate the levels associated to these
profiles, we must calculate the number of paths dominated by these profiles. According
to (a special case of) Kreweras' theorem, this amounts to calculating the determinants
of the matrices
$${t_{m-j+1} + 1 \choose 1 + j - i}_{+}$$
where  ${y \choose z}_{+}$ is defined as ${y \choose z}$ for $y \geq z \geq 0$
and $0$ otherwise. E.g., for
$m=m=10$, we easily calculate that the natural significance levels are
$5.41\times 10^{-6}$, $1.08\times 10^{-5}$, $2.71\times 10^{-5}$, 
$7.58\times 10^{-5}$, $2.27\times 10^{-4}$, 
 $7.14\times 10^{-4}$, $ 0.00232$, $0.00774$, $0.0263$, $0.0909$, and $0.318$.

{\bf Remark.}
In our theorem, we only considered the upper-tailed Smirnov test based on 
$D_{mn}^{+} = \sup_{t}\, (F_{m}(t) - G_{n}(t) )$. Of course, similar results exist for
the Smirnov tests based on $D_{mn}^{-} = \sup_{t}\, (G_{n}(t) - F_{m}(t) )$ or
$D_{mn} = \sup_{t}\, | F_{m}(t) - G_{n}(t) |$.

{\bf Acknowledgement.} We would like to thank John Einmahl for
stimulating discussions on this paper.

\end{document}